\documentclass{article}
\usepackage{algorithm}
\usepackage[noend]{algpseudocode}
\usepackage{amsmath}
\usepackage[numbers]{natbib}
\usepackage{graphicx} 

\usepackage{url}
\usepackage[hidelinks]{hyperref}
\usepackage[hyphenbreaks]{breakurl}

\usepackage{physics}

\title{Estimating the normal-inverse-Wishart distribution}
\author{Jonathan So, University of Cambridge \\ js2488@cam.ac.uk}
\date{May 2024}

\begin{document}

\maketitle

\section{Introduction}

The normal-inverse-Wishart (NIW) distribution is commonly used as a prior distribution for the mean and covariance parameters of a multivariate normal distribution. The family of NIW distributions is also a minimal exponential family. In this short note we describe a convergent procedure for converting from mean parameters to natural parameters in the NIW family, or---equivalently---for performing maximum likelihood estimation of the natural parameters given observed sufficient statistics. This is needed, for example, when using a NIW base family in expectation propagation \citep{minka2001family}.

\section{Normal-inverse-Wishart distribution}

The NIW distribution with parameters $\mu_0$, $\lambda$, $\Psi$ and $\nu$, has density \citep{wiki2023niw}
\begin{align}
    \text{NIW}(&\mu, \Sigma; \mu_0, \lambda, \Psi, \nu) = \nonumber\\
        &\frac{\lambda^{\frac{d}{2}}\lvert{\Psi}\rvert^{\frac{\nu}{2}}\lvert\Sigma\rvert^{-\frac{\nu + d + 2}{2}}}{(2\pi)^{\frac{d}{2}}2^{\frac{\nu d}{2}}\Gamma_d(\frac{\nu}{2})}\exp\left\{-\frac{1}{2}\text{Tr}\left(\Psi\Sigma^{-1}\right) - \frac{\lambda}{2}(\mu - \mu_0)\Sigma^{-1}(\mu - \mu_0)\right\} \label{eqn:pdf}.
\end{align}
We can recognise \eqref{eqn:pdf} as having the form of an exponential family distribution density, with statistic function
\begin{align}
    s\left(\mu, \Sigma\right) &= \left(-\frac{1}{2}\text{vec}\big(\Sigma^{-1}\big), \Sigma^{-1}\mu, -\frac{1}{2}\mu^\top\Sigma^{-1}\mu, -\frac{1}{2}\log\lvert\Sigma\rvert \right),
\end{align}
natural parameters
\begin{align}
    \eta &= \left(\text{vec}\big(\Psi + \lambda\mu_0\mu_0^\top\big), \lambda\mu_0, \lambda, \nu\right) \nonumber\\
        &=: (\eta_1, \eta_2, \eta_3, \eta_4), \label{eqn:natparams}
\end{align}
base measure
\begin{align}
    f(\mu, \Sigma) = \lvert\Sigma\rvert^{-\frac{d+2}{2}},
\end{align}
and log partition function
\begin{align}
    A(\eta) &= -\frac{d}{2}\log{\eta_3} - \frac{\eta_4}{2}\log\left\lvert\eta_1 - \frac{\eta_2\eta_2^\top}{\eta_3}\right\rvert + \frac{d}{2}\log(2\pi) + \frac{\eta_4 d}{2}\log 2 + \log \Gamma_d(\frac{\eta_4}{2}) \nonumber\\
    &= -\frac{d}{2}\log{\lambda} - \frac{\nu}{2}\log\left\lvert\Psi\right\rvert + \frac{d}{2}\log(2\pi) + \frac{\nu d}{2}\log 2 + \log \Gamma_d(\frac{\nu}{2}),
\end{align}
where we use $(a, b, \ldots)$ to denote concatenation of the vectors $a$, $b$, $\ldots$.

\section{Forward mapping}

We first derive the map from natural to mean parameters, a.k.a. the \emph{forward mapping}. Let $p$ be the density of the NIW distribution with natural parameters $\eta$; then, using the standard result for exponential families, $E_p[s(\mu, \Sigma)] = \nabla A(\eta)$ \citep{wainwright2008}, we have
\begin{align}
    \pdv{A}{\eta_1}
        &= \text{vec}\left\{-\frac{\eta_4}{2}\left(\eta_1 - \frac{\eta_2\eta_2^\top}{\eta_3}\right)^{-1}\right\} \nonumber\\
        &= -\frac{\nu}{2}\text{vec}(\Psi^{-1}) \nonumber\\
        &=: m_1, \\
    \pdv{A}{\eta_2} &= \frac{\eta_4}{\eta_3}\left(\eta_1 - \frac{\eta_2\eta_2^\top}{\eta_3}\right)^{-1}\eta_2 \nonumber\\
        &= \nu\Psi^{-1}\mu_0 \nonumber\\
        &=: m_2, \\
    \pdv{A}{\eta_3} &= -\frac{d}{2\eta_3} - \frac{\eta_4}{2\eta_3^2}\text{Tr}\left(\left(\eta_1 - \frac{\eta_2\eta_2^\top}{\eta_3}\right)^{-1}\eta_2\eta_2^\top\right) \nonumber\\
        &= -\frac{d}{2\lambda} - \frac{\nu}{2}\text{Tr}\big(\Psi^{-1}\mu_0\mu_0^\top\big) \nonumber\\
        &=: m_3, \\
    \pdv{A}{\eta_4} &= -\frac{1}{2}\log\left\lvert\eta_1 - \frac{\eta_2\eta_2^\top}{\eta_3}\right\rvert + \frac{d}{2}\log 2 + \frac{1}{2}\sum_{i=0}^{d-1}\psi\Big(\frac{\eta_4 - i}{2}\Big) \nonumber\\
        &= -\frac{1}{2}\log\left\lvert\Psi\right\rvert + \frac{d}{2}\log 2 + \frac{1}{2}\sum_{i=0}^{d-1}\psi\Big(\frac{\nu - i}{2}\Big) \nonumber\\
        &=: m_4,
\end{align}
so that $E_p[s(\mu, \Sigma)] = (m_1, m_2, m_3, m_4)$.

\section{Reverse mapping}

Let $M_1 = \text{vec}^{-1}(m_1)$; then, we have
\begin{align}
    \mu_0 &= (-2M_1)^{-1}m_2 \\
    \lambda &= -d(2m_3 + m_2^\top \mu_0)^{-1} \\
    \Psi &= \Big(-\frac{2}{\nu}M_1\Big)^{-1}.
\end{align}
Unfortunately $\nu$ cannot be solved for in closed form, but we can find it by performing a one-dimensional optimisation. We want to find $\nu$ such that
\begin{align}
    -\frac{1}{2}\log\left\lvert\Psi\right\rvert + \frac{d}{2}\log 2 + \frac{1}{2}\sum_{i=0}^{d-1}\psi\Big(\frac{\nu - i}{2}\Big) &= m_4,
\end{align}
or equivalently,
\begin{align}
    0 &= -\log\left\lvert\Psi\right\rvert + d\log 2 + \sum_{i=0}^{d-1}\psi\Big(\frac{\nu - i}{2}\Big) - 2m_4 \nonumber\\
    &= -\log\left\lvert\Big(-\frac{2}{\nu}M_1\Big)^{-1}\right\rvert + d\log 2 + \sum_{i=0}^{d-1}\psi\Big(\frac{\nu - i}{2}\Big) - 2m_4 \nonumber\\
    &= \log\left\lvert-2M_1\right\rvert - d\log\frac{\nu}{2} + \sum_{i=0}^{d-1}\psi\Big(\frac{\nu - i}{2}\Big) - 2m_4 \nonumber\\
    &=: f(\nu).
\end{align}
Taking first derivatives, we have
\begin{align}
    f'(\nu) &= -\frac{d}{\nu} + \frac{1}{2}\sum_{i=0}^{d-1}\psi^{(1)}\Big(\frac{\nu-i}{2}\Big).
\end{align}
Note that the density $p$ is only defined for $\nu > d-1$. Using the inequality \citep{wiki2023polygamma}
\begin{align}
    \frac{x + \frac{1}{2}}{x^2} \leq \psi^{(1)}(x),
\end{align}
we have
\begin{align}
    f'(\nu) &\geq -\frac{d}{\nu} + \sum_{i=0}^{d-1}\frac{(\nu - i) + 1}{(\nu - i)^2} \nonumber\\
        &\geq -\frac{d}{\nu} + \frac{d(\nu + 1)}{\nu^2} \nonumber\\
        &= \frac{d}{\nu^2} > 0,
\end{align}
hence $f'(\nu) > 0$, and $f(\nu)$ is strictly increasing for $\nu > d-1$. Taking second derivatives, we have
\begin{align}
    f''(\nu) &= \frac{d}{\nu^2} + \frac{1}{4}\sum_{i=0}^{d-1}\psi^{(2)}\Big(\frac{\nu-i}{2}\Big),
\end{align}
and using the inequality \citep{wiki2023polygamma}
\begin{align}
     \psi^{(2)}(x) \leq -\frac{1}{x^2} - \frac{1}{x^3},
\end{align}
we have the bound
\begin{align}
    f''(\nu) &\leq \frac{d}{\nu^2} - \frac{1}{4}\sum_{i=0}^{d-1}\frac{1}{\Big(\frac{\nu-i}{2}\Big)^2} + \frac{1}{\Big(\frac{\nu-i}{2}\Big)^3} \nonumber\\
        &= \frac{d}{\nu^2} - \sum_{i=0}^{d-1}\frac{1}{(\nu-i)^2} + \frac{2}{(\nu-i)^3} \nonumber\\
        &\leq \frac{d}{\nu^2} - \frac{d}{\nu^2} - \frac{2d}{\nu^3} \nonumber\\
        &= -\frac{2d}{\nu^3} < 0,
\end{align}
and so $f(\nu)$ is also strictly concave for $\nu > d-1$. The exponential family of NIW distributions is minimal, and so we know that there is a \emph{unique} $\nu$ for which $f(\nu) = 0$ \citep{wainwright2008}. Taken together, we have that $f(\nu)$ for $\nu > d-1$ is strictly increasing, strictly concave, and has a single root; this tells us that if we find a starting point $\nu$ for which $f(\nu) < 0$, the Newton-Raphson method is guaranteed to converge to the unique root. We can easily find such a point by taking any initial guess $\nu > d-1$ and halving the distance to between $\nu$ and $d-1$ until $f(\nu) < 0$. Once we have found $\mu_0$, $\lambda$, $\Psi$ and $\nu$, the natural parameters are given by \eqref{eqn:natparams}.

Pseudocode for the resulting procedure is given by Algorithm \ref{alg:niw_backwards_mapping}. $\nu_0$ is some initial guess for $\nu$ which could default to e.g. $\nu_0 = d$. $\epsilon$ is some small tolerance for detecting convergence.

\begin{algorithm}[ht]
    \caption{Normal-inverse-Wishart mean to natural parameter map}
    \label{alg:niw_backwards_mapping}
    \begin{algorithmic}
        \Require $(m_1, m_2, m_3, m_4)$
        \Require $\nu_0 > d-1$
        \Require $\epsilon > 0$
        \State $\nu \leftarrow \nu_0$
        \State $M_1 \leftarrow \text{vec}^{-1}(m_1)$
        \State $f(\nu) := \log\lvert -2M_1\rvert -d\log\frac{\nu}{2} + \sum_{i=0}^{d-1}\psi\Big(\frac{\nu-i}{2}\Big) - 2m_4$
        \While{$f(\nu) > 0$}
            \State $\nu \leftarrow (\nu - d + 1)/2$
        \EndWhile
        \While{$|f(\nu)| > \epsilon$}
            \State $\nu \leftarrow \nu - f(\nu)/f' (\nu)$
        \EndWhile
        \State $\mu_0 \leftarrow (-2M_1)^{-1}m_2$
        \State $\lambda \leftarrow -d(2m_3 + m_2^\top\mu_0)^{-1}$
        \State $\Psi \leftarrow (-\frac{2}{\nu}M_1)^{-1}$ \\
        \Return{$(\text{vec}(\Psi + \lambda\mu_0\mu_0^\top), \lambda\mu_0, \lambda, \nu)$}
    \end{algorithmic}
\end{algorithm}

\bibliography{main}

\end{document}